\documentclass[a4paper]{amsart}
\usepackage{graphicx} 
\usepackage{mathtools} 
\usepackage{cleveref} 
\usepackage[11pt]{moresize} 
\usepackage{appendix}

\usepackage[doi=false,url=false,maxnames=99,giveninits=true,isbn=false]{biblatex}
\usepackage{tikz} 
\tikzset{ext/.style={circle, draw,inner sep=1.4pt},int/.style={circle,draw,fill,inner sep=1.4pt},dot/.style={circle,fill,inner sep=0pt,outer sep=0pt}}
\usetikzlibrary{matrix} 

\newcommand{\K}{\mathbb{K}}
\newcommand{\sym}{\mathbb{S}}

\newcommand{\mV}{\mathrm{V}}
\newcommand{\mE}{\mathrm{E}}
\newcommand{\mB}{\mathrm{B}}

\newcommand{\mH}{\mathrm{H}}

\newcommand{\GC}{\mathrm{GC}}
\newcommand{\fGC}{\mathrm{fGC}}
\newcommand{\fGCc}{\mathrm{fGCc}}
\newcommand{\HGC}{\mathrm{HGC}}
\newcommand{\HoGC}{\mH_{\geq 0}\GC}

\newcommand{\grac}{\mathrm{grac}}

\newcommand{\coH}{\mathcal{H}}

\DeclareMathOperator{\sgn}{sgn}

\newtheorem{thm}{Theorem}
\newtheorem{prop}[thm]{Proposition}
\newtheorem{lem}[thm]{Lemma}
\newtheorem{cor}[thm]{Corollary}

\newtheorem{rem}[thm]{Remark}

\title{More on even graph cohomology cancellations}
\author{Andrija Tomorad}
\author{Marko Živković}
\date{November 2026}

\begin{document}

\begin{abstract}
    We provide an extra differential on even graph complexes whose differentials on odd pages of the corresponding spectral sequence are zero, implying that cancellations of even graph cohomology cancel only on even pages.
\end{abstract}

\maketitle

\section{Introduction}

\subsection{General remarks}

Graph complexes are graded vector spaces of formal linear combinations of isomorphism classes of some kind of graphs, with the standard differential defined by vertex splitting. The various graph cohomology theories are arguably some of the most fascinating objects in homological algebra. Each of graph complexes play a certain role in a subfield of homological algebra, algebraic topology or mathematical physics.

The simplest graph complex is introduced by Maxim Kontsevich in \cite{Kont1,Kont2}. It comes in versions $\GC_n$, where $n$ ranges over integers, and is spanned by simple connected graphs. The differential $\delta$ splits a vertex. See Section \ref{s:def} for the full definition. Some examples of graphs are:
$$
\begin{tikzpicture}[baseline=-1ex,scale=.7]
\node[int] (a) at (0,0) {};
\node[int] (b) at (1,0) {};
\draw (a) edge (b);
\draw (a) edge[bend left=40] (b);
\draw (a) edge[bend right=40] (b);
\end{tikzpicture}
\;,\quad
\begin{tikzpicture}[baseline=-1ex,scale=.6]
\node[int] (a) at (0,0) {};
\node[int] (b) at (90:1) {};
\node[int] (c) at (210:1) {};
\node[int] (d) at (330:1) {};
\draw (a) edge (b);
\draw (a) edge (c);
\draw (a) edge (d);
\draw (b) edge (c);
\draw (b) edge (d);
\draw (c) edge (d);
\end{tikzpicture}
\;.
$$

Since the differential increases both the number of vertices $v$ and the number of edges $e$ by $1$, the \emph{edge excess} $b\coloneq e-v$ is not changed\footnote{In literature an eqivalent value $e-v+1$ is often used, and called \emph{loop order} or \emph{Betti number}.}. Therefore, Kontsevich's graph complexes split into the product of sub-complexes with fixed $b$:
$$
(\GC_n,\delta)=\prod_{b\geq -1}(\mB_b\GC_n,\delta).
$$

Furthermore, the complexes $\GC_n$ and $\GC_{n'}$ are isomorphic up to some degree shifts if $n\equiv n' \text{ mod } 2$. Hence it suffices to understand two possible cases of $\GC_n$: \emph{even graph complex} for even $n$, and \emph{odd graph complex} for odd $n$.

The long standing open problem about graph complexes is to determine its cohomology. We know a little about it in general.

One of the standard techniques to gain some insights into the graph cohomology is to provide an extra differential. The technique has been introduced in \cite{DGC1} and later applied to another complexes or differentials in \cite{DGC2,DGC3,SED}. The rough idea is to twist the standard differential $\delta$ to another differential $\delta + \delta^{extra}$ that makes the complex (almost) acyclic. Then a spectral sequence is provided that on the first page has just the standard differential $\delta$, i.e.\ the cohomology of the first page of the spectral sequence is the cohomology $\coH(\GC_n,\delta)$ we are interested in. With the convergation of the spectral sequence properly checked, the acyclicity of the twisted complex implies that all classes of the standard cohomology cancel on further pages, i.e.\ they come in canceling pairs.

For Kontsevich's graph complexes an extra peace of differential $\delta^{extra}$ typically changes the edge excess $b=e-v$, and then the suitable spectral sequences is the one on $b$. So a class from $\mB_b\GC_n$ cancels another class from $\mB_{b'}\GC_n$ with $b'>b$.

For odd graph complex two extra differentials are known, one provided in \cite[Section 4]{DGC1} and another in \cite{SED}. The former one gives the acyclic complex
$$
(\GC_1, \delta + [m,\cdot])
$$
where the extra piece $[m,\cdot]$ is the Lie bracket with a Maurer-Cartan element $m$ that is a linear combination of graphs with odd edge excess. Details about Lie algebra structure on Kontsevich's graph complexes can be found e.g.\ in \cite{Will0}. Let us just note that the Lie bracket with a graph with odd edge excess changes the edge excess by even amount. Therefore the whole differential $\delta + [m,\cdot]$ cannot change the parity of edge excess. This implies that all cancellations come only on even pages of spectral sequence, i.e.\ a class from $\coH(\mB_b\GC_1,\delta)$ cancels another class from $\coH(\mB_{b'}\GC_1,\delta)$ with $b\equiv b' \text{ mod } 2$.

\subsection{Extra differentials for even graph cohomology}

In \cite[Section 3]{DGC1} an extra differential for even graph complex is provided, making the complex
\begin{equation}
\label{eq:nabla}
    (\GC_0, \delta + \nabla)
\end{equation}
acyclic. The extra peace $\nabla$ adds an edge in all possible ways. Here $\nabla$ increases edge excess by 1, giving no obvious reason that the cancellations appear only on certain pages. However, it has been immediately observed among known classes that the cancellations again appear only on even pages.

In \cite[Section 3.3.]{DGC3} the differential $\delta + \nabla$ on the version $\fGC_0$ that includes disconnected graphs has been conjugated to $\delta + D\nabla$, keeping the acyclicity. Here $D$ is \emph{deleting a vertex}, the map that deletes one vertex and reconnects edges that have been connected to that vertex to all other vertices in all possible ways. Note that $D\nabla$ increases edge excess by $2$, leading to the conclusion that all cancellations in the spectral sequence of
$$
(\fGC_0, \delta + D\nabla)
$$
come only on even pages. But there is a problem: deleting a vertex $D$ can both make connected graph disconnected and disconnected graph connected, and therefore this extra differential cannot be restricted to connected graph complex $\GC_0$. So it is possible that a connected class is actually cancelled by a disconnected one.

The latest achievement towards restricting odd-page cancellations in even graph cohomology $\coH(\GC_0,\delta)$ is done by Willwacher in \cite[Corollary 4.4.]{TCon}, where he proved that differentials on pages 1, 3 and 5 are zero by checking that a connected class cannot be mapped to a disconnected class up to page 5.

The purpose of this paper is to provide an extra differential for even graph cohomology with no odd-page cancellations altogether.

\subsection{Hairy graph complexes}

To do this we use techniques similar to \cite{SED}, i.e.\ instead of an extra differential directly on $\GC_0$ we provide an extra differential on its quasi-isomorphic picture consisting of hairy graphs.

Hairy graph complexes $\HGC_{m,n}$, where $m$ and $n$ range over integers, are spanned by connected graphs with external legs, also called ``hairs'' (see Section \ref{s:def} for the definition). At least one hair is required. Similarly to Kontsevich's graph complexes, there are essentially 4 types depending on the parities of $m$ and $n$.
Some examples of hairy graphs are:
\begin{equation}
\label{eq:HGexample}
\begin{tikzpicture}[baseline=-1ex,scale=.7]
\node[int] (a) at (0,0) {};
\draw (a) edge (90:.3);
\draw (a) edge (120+90:.3);
\draw (a) edge (240+90:.3);
\end{tikzpicture}
\;,\quad
\begin{tikzpicture}[baseline=-2ex,scale=.6]
\node[int] (v1) at (-1,0){};
\node[int] (v2) at (0,.8){};
\node[int] (v3) at (1,0){};
\node[int] (v4) at (0,-.8){};
\draw (v1)  edge (v2);
\draw (v1)  edge (v4);
\draw (v2)  edge (v3);
\draw (v3)  edge (v4);
\draw (v4)  edge (v2);
\draw (v1)  edge (-1.3,0);
\draw (v3)  edge (1.3,0);
\end{tikzpicture}
\;,\quad
\begin{tikzpicture}[baseline=-1ex,scale=.6]
\node[int] (a) at (0,0) {};
\node[int] (b) at (90:1) {};
\node[int] (c) at (210:1) {};
\node[int] (d) at (330:1) {};
\draw (b) edge (0,1.3);
\draw (a) edge (b);
\draw (a) edge (c);
\draw (a) edge (d);
\draw (b) edge (c);
\draw (b) edge (d);
\draw (c) edge (d);
\end{tikzpicture}
\;.
\end{equation}

There are multiple extra differentials on hairy graph complexes. One is provided in \cite{TW1} on $\HGC_{n-1,n}$ to be Lie-bracket with a Maurer-Cartan element $\sigma$ \eqref{eq:sigma}.
Hairy graph complexes are equipped with a distinguished Lie algebra structure, see the definition in Subsection \ref{ss:Lie}.
The resulting complex is (almost) quasi-isomorphic to the Kontsevich's graph complex, i.e.\ there is a quasi-isomorphism
\begin{equation}
\label{eq:TWqi}
    \K \oplus (\GC_{n},\delta) \to \left(\HGC_{n-1,n}, \delta+[\sigma,\cdot]\right).
\end{equation}

Another extra differential, this time on $\HGC_{-1,n}$, is $\Delta$ that connects a hair into an edge. Here we know that
\begin{equation}
\label{eq:Delta}
    (\HGC_{-1,n},\delta+\Delta)
\end{equation}
is acyclic (\cite{DGC3}, conjectured in \cite{DGC2}).

\subsection{Sketch of the proof}

The main idea of our result is extending \eqref{eq:TWqi} to the explicit quasi-isomorphism
\begin{equation}
    \K \oplus (\GC_{0},\delta)\oplus (\GC_{0},\delta) \to \left(\mH_{\geq 0}\GC_{-1,0}, \delta+[\sigma,\cdot]\right),
\end{equation}
where $\mH_{\geq 0}\GC_{-1,0}=\GC_0[-1]\oplus\HGC_{-1,0}$ is the hairy graph complex with hair-less graphs also allowed. We show that the quasi-isomorphism behaves well (commutes) with further adding extra differentials
\begin{equation}
    \K \oplus (\GC_{0},\delta+\nabla)\oplus (\GC_{0},\delta-\nabla) \to \left(\mH_{\geq 0}\GC_{-1,0}, \delta+[\sigma,\cdot]+\Delta\right).
\end{equation}
It is interesting that on the the left-hand side we have twice our already known extra differential $\nabla$ from \eqref{eq:nabla}. So the left-hand side is acyclic, up to one class provided by $\K$, and therefore is the right-hand side too.

Now, on the right-hand side we construct another spectral sequence such that the first differential is $\delta+\Delta$. If the number of hairs is denoted by $h$, this spectral sequence is on the number $f\coloneq e+h-v$, called \emph{edge-hair excess}.

In this spectral sequence we have an extra differential $[\sigma,\cdot]$ that changes edge-hair excess always by an even amount. Since the whole complex is acyclic, all classes of $\coH\left(\mH_{\geq 0}\GC_{-1,0},\delta+\Delta\right)$ eventually cancel, and all cancellations are between classes of the same parity of edge-hair excess.

The acyclicity of \eqref{eq:Delta} implies that the inclusion
\begin{equation}
    (\GC_0,\delta)\hookrightarrow(\mH_{\geq 0}\GC_{-1,0},\delta+\Delta)
\end{equation}
is quasi-isomorphism, so the cancellations in the spectral sequence are indeed cancellations of classes from even graph cohomology $(\GC_0,\delta)$.
Note that the inclusion maps a sub-complex with fixed edge excess $b$ into the sub-complex with edge-hair excess $b$, so the cancellations in $(\GC_0,\delta)$ are between classes of the same parity of edge excess.

\subsection*{Structure of the paper}

Section \ref{s:def} recalls needed definitions and results from the literature. In Section \ref{s:hairs} we discuss the extra differential $[\sigma,\cdot]$ and its implications, while in Section \ref{s:extra} we discuss the extra differentials $\nabla$ and $\Delta$ and deduce the final result. The Appendix \ref{a:pict} recalls what we know about the even graph cohomology by now, and the Appendix \ref{a:ex} gives an example of canceling classes using our new result.

\subsection*{Declaration of AI usage}

The AI tools have not been used for the mathematical content in this paper.
The free version of OpenAI’s ChatGPT has been used for references and LaTeX formatting.

\subsection*{Acknowledgments}

This work was supported in part by the Croatian Science Foundation under the project number HRZZ-IP-2026-2389 (TEPiC), and in part by the European Union - NextGenerationEU through the National Recovery and Resilience Plan 2021–2026, via an institutional grant of the University of Zagreb Faculty of Science, IK IA 1.1.3, Impact4Math.

\section{Graph complexes}
\label{s:def}

In this section we recall definitions of graph complexes and some other notions needed in this paper, and fix our notation.

Standard Kontsevich's graph complexes $\GC_n$ and hairy graph complexes $\HGC_{m,n}$ are in general defined e.g.\ in \cite{DGC3}.
In this paper we are only interested in the case when $m=-1$ and $n=0$, that is $\GC_0$ and $\HGC_{-1,0}$. For simplicity, we will use the shorter notations $\GC:=\GC_0$ and $\HGC:=\HGC_{-1,0}$.

We work over a field $\K$ of characteristic zero. All vector spaces and differential graded vector spaces are assumed to be $\K$-vector spaces. 

\subsection{From graphs to space of co-invariants}
Consider the set $\bar\mV_v\bar\mE_e\bar\mH_h\grac$ containing graphs that:
\begin{itemize}
\item are connected;
\item have $v>0$ distinguishable vertices;
\item have $e\geq 0$ distinguishable undirected edges between different vertices\footnote{Thus, edges that start and end at the same vertex, often called ``tadpoles'', are not allowed.};
\item have $h\geq 0$ distinguishable hairs attached to some vertices;
\end{itemize}
For some pictures of such graphs see \eqref{eq:HGexample}.

Let
\begin{equation}
\bar\mV_v\bar\mE_e\bar\mH_h\GC:=\langle\bar\mV_v\bar\mE_e\bar\mH_h\grac\rangle[-e-1]
\end{equation}
be the vector space of formal series of $\bar\mV_v\bar\mE_e\bar\mH_h\grac$ with coefficients in $\K$. It is a graded vector space with a non-zero term only in degree $d=e+1$.

There is a natural right action of the group $\sym_v\times \sym_h\times \sym_e$ on $\bar\mV_v\bar\mE_e\bar\mH_h\grac$, where $S_v$ permutes vertices, $S_h$ permutes hairs and $S_e$ permutes edges.
Let $\sgn_e$ be one-dimensional representations of $S_e$, where the odd permutation reverses the sign. It can be considered as a representation of the whole product $\sym_v\times \sym_h\times \sym_e$. Let us consider the space of co-invariants:
\begin{equation}
\mV_v\mE_e\mH_s\GC:=
\left(\bar\mV_v\bar\mE_e\bar\mH_h\GC\otimes\sgn_e\right)_{\sym_v\times \sym_h\times \sym_e}
\end{equation}
The operation of taking co-invariants effectively removes numbering on vertices and hairs, and removes numbering of edges up to sign.

\begin{rem}
Sign changes induced by switching edges imply that in $\mV_v\mE_e\mH_h\GC$ there are no multiple edges, i.e.\ more than one edges that start and end on same pair of vertices.
\end{rem}

\subsection{The differential}
\label{ss:dif}
For $\Gamma\in\bar\mV_v\bar\mE_e\bar\mH_h\GC$ let ``splitting a vertex'' be
\begin{equation}
    \sigma(\Gamma)\coloneq\sum_{x\in V(\Gamma)}\sigma_x(\Gamma),
\end{equation}
where $V(\Gamma)$ is the set of vertices of $\Gamma$, and $\sigma_x$ is inserting
\begin{tikzpicture}[scale=.5]
 \node[int] (a) at (0,0) {};
 \node[int] (b) at (1,0) {};
 \draw (a) edge (b);
\end{tikzpicture}
instead of the vertex $x$ and summing over all possible ways of connecting the edges that have been connected to $x$ to the new two vertices.
Let ``adding an antenna'' be
\begin{equation}
    \alpha(\Gamma)\coloneq\sum_{x\in V(\Gamma)}\alpha_x(\Gamma),
\end{equation}
where $\alpha_x$ adds an edge on $x$ towards new 1-valent vertex.
Let ``extracting a hair'' be
\begin{equation}
    \epsilon(\Gamma)\coloneq\sum_{k\in H(\Gamma)}\epsilon_k(\Gamma),
\end{equation}
where $H(\Gamma)$ is the set of hairs of $\Gamma$, and $\epsilon_k$ switches the hair with an edge towards a new vertex with that hair.

One checks that the maps $\sigma$, $\alpha$ and $\epsilon$ can be extended to the space of co-invariants $\mV_v\mE_e\mH_h\GC$. To precisely define the sign we set that the new edge gets the last label in a representative.

The differential on $\mV_v\mE_e\mH_h\GC$ is defined as:
\begin{equation}
 \delta(\Gamma)\coloneq\frac{1}{2}\sigma(\Gamma)-\alpha(\Gamma)-\epsilon(\Gamma).
\end{equation}
Note that, unless the whole graph has a single vertex, $\alpha$ and $\epsilon$ cancel terms from $\sigma$ where all edges end hairs, apart from possibly one hair, are connected to the same new vertex.

\subsection{Complexes}

The differential does not change the number of hairs, so for each $h\geq 0$ we can define a graph complex
$$
\left(\mH_h\GC,\delta\right)
$$
where
$$
\mH_h\GC\coloneq\prod_{v\geq 1, e\geq 0}\mV_v\mE_e\mH_h\GC.
$$

Also the complex does not change \emph{edge excess} $b\coloneq e-v$ so those complexes can be split as the direct product of sub-complexes
$$
(\mH_h\GC,\delta)=\left(\prod_{b\geq -1}\mH_h\mB_b\GC,\delta\right)
$$
where
$$
\mH_h\mB_b\GC\coloneq\prod_{v\geq 1}\mV_v\mE_{b+v}\mH_h\GC.
$$

\subsection{Kontsevich's graph complex}
\label{ss:comp}

Kontsevich's connected full graph complex is
\begin{equation}
(\fGCc,\delta)=(\mH_0\GC[1],\delta).
\end{equation}

It is common to define sub-complexes of $\fGCc$ spanned by graphs whose vertices are at least $2$- or $3$-valent, denoted by $\fGCc^2$ and $\fGCc^3$ respectively. Note that the differential cannot ruin those valence conditions, so they are indeed sub-complexes.

It is a well known result (see e.g.\ \cite{Will0}) that
\begin{equation}
    \coH(\fGCc,\delta)=\coH(\fGCc^2,\delta)=\coH(\mB_0\fGCc^2,\delta)\oplus\coH(\fGCc^3,\delta),
\end{equation}
where $\coH(\mB_0\fGCc^2,\delta)$ consists of loop classes represented by loops $L_{4k+1}$ for positive integers $k$. $L_{4k+1}$ is a graph with $4k+1$ vertices and $4k+1$ edges arranged in one loop, e.g.\ \eqref{eq:L5}. Also $(\mB_{-1}\fGCc^2,\delta)$, spanned by trees, is acyclic.

From now on we abbreviate\footnote{In literature it is more common to abbreviate $\GC=\fGCc^3$. However, for our purposes it is simpler to use $\fGCc^2$ because the cancellations happen in $\coH(\fGCc,\delta)=\coH(\fGCc^2,\delta)$.
}
\begin{equation}
    \GC\coloneq\fGCc^2.
\end{equation}

\subsection{Hairy graph complex}

The hairy graph complex is
\begin{equation}
(\HGC,\delta)=\left(\prod_{h\geq 1}\mH_h\GC,\delta\right).
\end{equation}
We also need a combined complex
\begin{equation}
\left(\mH_{\geq 0}\GC,\delta\right)=\left(\prod_{h\geq 0}\mH_h\GC,\delta\right)=
(\fGCc[-1],\delta)\oplus(\HGC,\delta).
\end{equation}

It is known that the tree part $(\mB_{-1}\HGC,\delta)$ is $1$-dimensional, spanned by $\Sigma_3$ (e.g.\ \cite[Subsection 3.5]{AT} or \cite{DGC3}).

In the rest of the article we will often define functions on a particular vertex or a hair of a graph, and often take sums of them. Technically it can only be done on graphs $\Gamma\in \bar\mV_v\bar\mE_e\bar\mH_h\GC$, i.e.\ before taking co-invariants. However, by abuse of notation, we say we do that on $\Gamma\in\HoGC$ (or another graph complex), meaning that we do it on a representative of $\Gamma$ in the space of co-invariants.
When there are sums over all vertices or hairs, one can check that the definition can always be extended to the space of co-invariants.

\subsection{Lie algebra structure on the hairy graph complex}
\label{ss:Lie}

Given two hairy graphs $\Gamma, \Gamma'\in\HoGC$ of homogeneous degree, we define the Lie bracket as
\begin{equation}
    [\Gamma,\Gamma']=\Gamma\circ\Gamma'-(-1)^{|\Gamma||\Gamma'|}\Gamma'\circ\Gamma
\end{equation}
where $\circ$ is the operation of connecting one hair of $\Gamma$ to vertices of $\Gamma'$ in all possible ways. 
To precisely define the sign we set that the representative of $\Gamma\circ\Gamma'$ has edges of $\Gamma'$ labeled first, in their original order, then the connecting edge gets the label, and finally the edges of $\Gamma$ are labeled in their order. For more details see eg.\ \cite{TW1}.

\subsection{Shifting the differential}

For an integer $a\geq 1$ the \emph{corolla} $\Sigma_a$ is a graph with a single vertex and $a$ hairs on it:
\begin{equation}
\Sigma_a=
\begin{tikzpicture}[baseline=-.5ex]
\node[int] (a) at (0,0) {};
\draw (a) edge (90:.2);
\draw (a) edge (360/7+90:.2);
\node[dot] at (2*360/7+90:.15) {};
\node[dot] at (2.5*360/7+90:.15) {};
\node[dot] at (3*360/7+90:.15) {};
\draw (a) edge (4*360/7+90:.2);
\draw (a) edge (5*360/7+90:.2);
\draw (a) edge (6*360/7+90:.2);
\end{tikzpicture}.
\end{equation}
Observe that for all $\Gamma\in\HoGC$
\begin{equation}
\label{eq:AlphaEpsilonSigma}
    (\alpha+\epsilon)(\Gamma)=[\Sigma_1,\Gamma].
\end{equation}
Let
\begin{equation}
\label{eq:sigma}
    \Sigma:=\sum_{k=1}^\infty\frac{1}{(2k+1)!}\Sigma_{2k+1}.
\end{equation}
One can check that $\Sigma$ is a Maurer-Cartan element in dg Lie algebra $(\HoGC,\delta)$, i.e.\ it satisfies:
\begin{equation*}
    -\delta\Sigma+\frac{1}{2}[\Sigma,\Sigma]=0.
\end{equation*}
It can be used to shift the differential to the new differential
\begin{equation}
    \delta_\Sigma\coloneq\delta - [\Sigma,\cdot].
\end{equation}
Using \eqref{eq:AlphaEpsilonSigma} we get
\begin{equation}
    \delta_\Sigma=\sigma-[\Sigma',\cdot].
\end{equation}
where
\begin{equation}
    \Sigma' \coloneq \Sigma_1+\Sigma = \sum_{k=0}^\infty\frac{1}{(2k+1)!}\Sigma_{2k+1}.
\end{equation}

\section{Adding hairs and hairy structures}
\label{s:hairs}

\subsection{Decorated graphs}
\label{ss:decor}

On $\HoGC$ we define the map $\chi_x$ that adds a hair on vertex $x$. For $\Gamma\in\mH_0\GC$ we can exponentiate it in the standard way:
\begin{equation}
    e^{\pm\chi_x}(\Gamma)\coloneq\sum_{i=0}^\infty\frac{(\pm 1)^i}{i!}\chi_x^i(\Gamma).
\end{equation}

A \emph{Decorated graph} is a pair $(\Gamma,s)$ where $\Gamma\in\mH_0\GC$ is a graph and $s\colon V(\Gamma)\rightarrow\{+1,-1\}$, $x\mapsto s_x$. Decorated graph $(\Gamma,s)$ represents the graph
\begin{equation}
    \left(\prod_{x\in V(\Gamma)}e^{s_x\chi_x}\right)(\Gamma)\in\HoGC.
\end{equation}
It is a graph on whose all vertices we acted by either $e^{\chi_x}$ or $e^{-\chi_x}$. Decorated graph $(\Gamma,s)$ can be understood as a hairless graph with vertices decorated by $s_x$, an example being:
\begin{equation*}
\begin{tikzpicture}[baseline=-1ex,scale=.8]
\node[int, label={[xshift=4pt,yshift=-4pt]:$+$}] (a) at (0,0) {};
\node[int, label=above:$-$] (b) at (90:1) {};
\node[int, label=below left:$+$] (c) at (210:1) {};
\node[int, label=below right:$-$] (d) at (330:1) {};
\draw (a) edge (b);
\draw (a) edge (c);
\draw (a) edge (d);
\draw (b) edge (c);
\draw (b) edge (d);
\draw (c) edge (d);
\end{tikzpicture}
\end{equation*}

\subsection{The differential on decorated graphs}

\begin{lem}
\label{lem:split}
    For a decorated graph $(\Gamma,s)\in\HoGC$ the splitting a vertex acts as
    $$
    \sigma_x(\Gamma,s)=(\sigma_x(\Gamma),\sigma_x(s))
    $$
    where $\sigma_x(s)$ is equal to $s$ on all vertices $y\neq x$, and is equal to $s_x$ on both new vertices that come from splitting $x$.
\end{lem}
\begin{proof}
    Let us call the newly formed vertices $x$ and $x'$. Then
    \begin{equation*}
        \begin{split}
            \sigma_x(\Gamma,s)
            &= \sigma_x\left(\left(\prod_{y\in V(\Gamma)}e^{s_y\chi_y}\right)(\Gamma)\right) \\
            &= \left(\prod_{y\in V(\Gamma)\setminus\{x\}}e^{s_y\chi_y}\right)
            \sigma_x\left(e^{s_x\chi_x}(\Gamma)\right) \\
            &= \left(\prod_{y\in V(\Gamma)\setminus\{x\}}e^{s_y\chi_y}\right)
            \sigma_x\left(\sum_{i=0}^\infty\frac{s_x^i}{i!}\chi_x^i(\Gamma)\right) \\
            &= \left(\prod_{y\in V(\Gamma)\setminus\{x\}}e^{s_y\chi_y}\right)
            \sum_{i=0}^\infty\frac{s_x^i}{i!}\sigma_x\left(\chi_x^i(\Gamma)\right) \\
            &= \left(\prod_{y\in V(\Gamma)\setminus\{x\}}e^{s_y\chi_y}\right)
            \sum_{i=0}^\infty\frac{s_x^i}{i!}\sum_{j=0}^i\binom{i}{j}\chi_x^j\chi_{x'}^{i-j}\sigma_x(\Gamma) \\
            &= \left(\prod_{y\in V(\Gamma)\setminus\{x\}}e^{s_y\chi_y}\right)
            \sum_{j=0}^\infty\frac{s_x^j}{j!}\sum_{k=0}^\infty\frac{s_x^k}{k!}\chi_x^j\chi_{x'}^k\sigma_x(\Gamma) \\
            &= \left(\prod_{y\in V(\Gamma)\setminus\{x\}}e^{s_y\chi_y}\right)
            e^{s_x\chi_x}e^{s_x\chi_{x'}}\sigma_x(\Gamma) \\
            &= (\sigma_x(\Gamma),\sigma_x(s)).
        \end{split}
    \end{equation*}
\end{proof}

Let $\Gamma\in\HoGC$, $x\in V(\Gamma)$ and $i\geq 0$ an integer. Define
\begin{equation}
    \phi^i_x(\Gamma)\coloneq\chi_{x'}^i\alpha_x(\Gamma)
\end{equation}
where $x'$ is the end of ``antenna'' added by $\alpha$. The map $\phi^i_x(\Gamma)$ can be called ``adding a flower'' and adds an edge from $x$ to the new vertex, and adds $i$ hairs on the new vertex.

It is easily seen that for $\Gamma\in\mH_0\GC$
\begin{equation*}
    \chi_x^i(\Gamma)\circ\Sigma_i = -(-1)^{|\Gamma|}i\phi_x^i(\chi^{i-1}(\Gamma)),
\end{equation*}
implying that
\begin{equation}
    e^{\pm\chi_x}(\Gamma)\circ\Sigma_i=\mp(-1)^{|\Gamma|}\phi_x^i(e^{\pm\chi_x}(\Gamma)).
\end{equation}
Adding hairs to different vertices commute, so for a decorated graph $(\Gamma,s)$
\begin{equation}
\label{eq:SigmaPhi2}
    \begin{split}
        (\Gamma,s)\circ\Sigma_i
        &= \left(\prod_{y\in V(\Gamma)}e^{s_y\chi_y}\right)(\Gamma)\circ\Sigma_i \\
        &= -(-1)^{|\Gamma|}\sum_{x\in V(\Gamma)}s_x\phi_x^i\left(\left(\prod_{y\in V(\Gamma)}e^{s_y\chi_y}\right)(\Gamma)\right) \\
        &= -(-1)^{|\Gamma|}\sum_{x\in V(\Gamma)}s_x\phi_x^i(\Gamma,s).
    \end{split}
\end{equation}

Also let
\begin{equation}
    \phi^i(\Gamma)\coloneq\sum_{x\in V(\Gamma)}\phi^i_x(\Gamma).
\end{equation}
Then for all $\Gamma\in\HoGC$ it holds that
\begin{equation}
\label{eq:SigmaPhi1}
    \Sigma_i\circ\Gamma = i\phi^{i-1}(\Gamma).
\end{equation}

\begin{lem}
\label{lem:shift}
    Let $(\Gamma,s)$ be a decorated graph. Then
    $$
    [\Sigma',(\Gamma,s)] = \sum_{x\in V(\Gamma)}(\alpha_x(\Gamma),\sigma_x(s)).
    $$
    where $\sigma_x(s)$ is as in Lemma \ref{lem:split}.
\end{lem}
\begin{proof}
    \begin{equation*}
        \begin{split}
            [\Sigma',(\Gamma,s)]
            &= \sum_{k=0}^\infty\frac{1}{(2k+1)!}[\Sigma_{2k+1},(\Gamma,s)] \\
            &= \sum_{k=0}^\infty\frac{1}{(2k+1)!}\left(\Sigma_{2k+1}\circ(\Gamma,s)-(-1)^{|\Gamma|}(\Gamma,s)\circ\Sigma_{2k+1}\right) \\
            &\overset{\eqref{eq:SigmaPhi1}\eqref{eq:SigmaPhi2}}{=}
            \sum_{k=0}^\infty\frac{1}{(2k+1)!}
            \left((2k+1)\phi^{2k}(\Gamma,s)
            +\sum_{x\in V(\Gamma)}s_x\phi_x^{2k+1}(\Gamma,s)\right) \\
            &= \sum_{x\in V(\Gamma,s)}\left(\sum_{k=0}^\infty\frac{1}{(2k)!}\phi_x^{2k}(\Gamma,s)
            \pm \sum_{k=0}^\infty\frac{s_x}{(2k+1)!}\phi^{2k+1}(\Gamma,s)\right) \\
            &= \sum_{x\in V(\Gamma)}\sum_{i=0}^\infty\frac{s_x^i}{i!}\phi_x^i(\Gamma,s) \\
            &= \sum_{x\in V(\Gamma)}(\alpha(\Gamma),\sigma_x(s)).
        \end{split}
    \end{equation*}
\end{proof}

\subsection{Decorated complex}

Let
\begin{equation}
    \GC^\pm\coloneq\langle\{(\Gamma,s)\mid\Gamma\in\GC[-1]\}\rangle\subseteq\HoGC.
\end{equation}
Note that in the space $\GC^\pm$ underlying graphs have at least 2-valent vertices.

\begin{prop}
\label{prop:GCpm}
    On a decorated graph $(\Gamma,s)\in\GC^\pm$ the differential $\delta_\Sigma$ acts as
    $$
    \delta_\Sigma(\Gamma,s) = \sum_{x\in V(\Gamma)}(\delta_x(\Gamma),\sigma_x(s))
    $$
    where $\delta_x=\frac{1}{2}\sigma_x-\alpha_x$ and $\sigma_x(s)$ is as in Lemma \ref{lem:split}.

    Therefore $\GC^\pm$ forms a sub-complex of $(\HoGC,\delta_\Sigma)$.
\end{prop}
\begin{proof}
    Using Lemmas \ref{lem:split} and \ref{lem:shift} it is
    \begin{equation*}
        \begin{split}
            \delta_\Sigma(\Gamma,s)
            &= \frac{1}{2}\sigma(\Gamma,s)-[\Sigma',(\Gamma,s)] \\
            &= \frac{1}{2}\sum_{x\in V(\Gamma)}(\sigma_x(\Gamma),\sigma_x(s))
            -\sum_{x\in V(\Gamma)}(\alpha_x(\Gamma),\sigma_x(s)) \\
            &= \sum_{x\in V(\Gamma)}(\delta_x(\Gamma),\sigma_x(s)).
        \end{split}
    \end{equation*}
\end{proof}

\subsection{Uniform adding hairs}

On $\HoGC$ we define degree zero map:
\begin{equation}
\chi(\Gamma)\coloneq\sum_{x\in V(\Gamma)}\chi_x(\Gamma).
\end{equation}
 For $\Gamma\in\mH_0\GC$ we have
\begin{equation}
    e^{\pm\chi}(\Gamma)
    =\left(\prod_{x\in V(\Gamma)}e^{\pm\chi_x}\right)(\Gamma)
    =(\Gamma,s_\pm)\in\GC^\pm,
\end{equation}
where $s_+$, respectively $s_-$, is $+1$, respectively $-1$, on all vertices of $\Gamma$.

Let $\GC^+$ be spanned by all decorated graphs $(\Gamma,s_+)$ for $\Gamma\in\GC$, i.e.\ by graphs in $\GC^\pm$ with all vertices decorated with $+$. Proposition \ref{prop:GCpm} implies that $\GC^+$ forms a sub-complex of $\GC^\pm$ and that
\begin{equation}
\label{eq:qi+}
    e^\chi\colon(\GC,\delta)\rightarrow(\GC^+,\delta_\Sigma)
\end{equation}
is an isomorphism.
Similarly we define sub-complex $\GC^-\subseteq\GC^\pm$ and we have an isomorphism
\begin{equation}
\label{eq:qi-}
    e^{-\chi}\colon(\GC,\delta)\rightarrow(\GC^-,\delta_\Sigma).
\end{equation}

\subsection{Revisiting the proof of \cite[Theorem 3]{TW1}}

Composed with inclusions, those isomorphisms give us two maps from $(\GC,\delta)\rightarrow(\HoGC,\delta_\Sigma)$. They can be combined into one map $(\GC,\delta)\oplus(\GC,\delta)\rightarrow(\HoGC,\delta_\Sigma)$:
\begin{equation}
\label{eq:map}
(\Gamma_1,\Gamma_2)\mapsto e^\chi(\Gamma_1)+e^{-\chi}(\Gamma_2).
\end{equation}

Recall that in \cite[Theorem 3]{TW1} it has been proven that there exists a quasi-isomorphism $(\GC,\delta)\rightarrow(\HGC,\delta_\Sigma)$. However, the proof itself first shows it for the extended situation $(\GC,\delta)\oplus(\GC,\delta)\rightarrow(\HoGC,\delta_\Sigma)$. In this subsection we go deeply into this proof to check that the needed quasi-isomorphism is indeed \eqref{eq:map}. For the sake of completeness we will recall all the needed original arguments.

\begin{thm}[{Partially contained in the proof of \cite[Theorem 3]{TW1}}]
\label{thm:qi}
    The map $\K\oplus(\GC,\delta)\oplus(\GC,\delta)\rightarrow(\HoGC,\delta_\Sigma)$ defined with
    \begin{equation*}
        (\lambda,\Gamma_1,\Gamma_2)\mapsto \lambda\Sigma_3+e^\chi(\Gamma_1)+e^{-\chi}(\Gamma_2)
    \end{equation*}
    is a quasi-isomorphism.
\end{thm}
\begin{proof}
Note that neither the differential $\delta_\Sigma$ nor the morphisms $e^{\pm\chi}$ change edge access, so we can separate the trees (edge access part $-1$). We know that $(\mB_{-1}\GC,\delta)$ is acyclic and that $\coH(\mB_{-1}\HoGC,\delta)$ is spanned by $\Sigma_3$, c.f.\ Subsection \ref{ss:comp}. Easy spectral sequence argument, left to the reader, implies that $\coH(\mB_{-1}\HoGC,\delta_\Sigma)$ is also spanned by $\Sigma_3$. So to prove the Theorem it is enough to prove that $(\GC,\delta)\oplus(\GC,\delta)\rightarrow(\mB_{\geq 0}\HoGC,\delta_\Sigma)$ defined with
    \begin{equation*}
        (\Gamma_1,\Gamma_2)\mapsto e^\chi(\Gamma_1)+e^{-\chi}(\Gamma_2)
    \end{equation*}
is quasi-isomorphism. Because of isomorphisms \eqref{eq:qi+} and \eqref{eq:qi-} it is enough to prove that the inclusion
$$
\GC^+\oplus\GC^-\hookrightarrow\mB_{\geq 0}\HoGC
$$
is a quasi-isomorphism. This follows from the following Propositions \ref{prop:1} and \ref{prop:2}.
\end{proof}

We use the following well known lemma.

\begin{lem} \label{lem:q-isom}
    Let $(V,\delta)$ be a (co-)chain complex with a descending bounded filtration. Let $(V,\delta_0)$ be its associated graded complex. If $(W,\delta)$ is a sub-complex of $(V,\delta)$ such that the inclusion induces isomorphism of the associated graded homologies $\coH(W,\delta_0)\rightarrow\coH(V,\delta_0)$, then the inclusion of the full complexes $(W,\delta)\hookrightarrow(V,\delta)$ is a quasi-isomorphism.
\end{lem}
\begin{proof}
    On the mapping cone of the inclusion $(W,\delta)\hookrightarrow(V,\delta)$ we consider the same filtration and the associated spectral sequence. The associated graded complex is exactly the mapping cone of the inclusion $(W,\delta_0)\hookrightarrow(V,\delta_0)$.

    If the induced map $\coH(W,\delta_0)\rightarrow\coH(V,\delta_0)$ is an isomorphism, the mapping cone of $(W,\delta_0)\hookrightarrow(V,\delta_0)$ is acyclic. For the bounded filtration the spectral sequence converges correctly, so the whole mapping cone of $(W,\delta)\hookrightarrow(V,\delta)$ is acyclic, and hence the inclusion is a quasi-isomorphism.
\end{proof}

\begin{prop}
\label{prop:1}
    The inclusion
    $$
    \GC^\pm\hookrightarrow\mB_{\geq 0}\HoGC
    $$
    is a quasi-isomorphism.
\end{prop}
\begin{proof}
    We define the \emph{core} of the graph $\Gamma$ as its smallest possible sub-graph which has the same loop order as $\Gamma$. That is, it is the graph obtained by removing all hairs, and recursively removing all the univalent vertices. Each graph can be perceived as a compound of its core and hairy trees decorating its vertices. The idea is that $\GC^\pm$ replaces the decorating trees with simpler structures which have the same contribution to the co-homology.

Now we consider the filtration by number of vertices in the core on $\mB_{\ge0}\mH_{\ge 0}\GC$. In that case, the first differential on the filtered complex is $\delta_\Sigma^0$ being the part of $\delta_\Sigma$ that does not change the core. Therefore, the associated graded complex splits as:
$$
 \left(B_{\ge 0}H_{\ge 0}\GC,\delta_\Sigma^0\right) = \prod_{\Gamma} \left(V_\Gamma,\delta_\Sigma^0\right)
$$
where for a fixed core $\Gamma$, $V_\Gamma$ is the complex spanned by all graphs with that core. We do the same filtration and splitting for the sub-complex $\GC^\pm$:
$$
\left(\GC^\pm,0\right) = \prod_{\Gamma} \left(V'_\Gamma,0\right).
$$
Because of Proposition \ref{prop:GCpm} the differential on $\GC^\pm$ must change the core, so the first differential here is zero. Note that the core of a decorated graph $(\Gamma,s)\in\GC^\pm$ is exactly $\Gamma$, so $V'_\Gamma$ is the complex spanned by $\{(\Gamma,s)\mid s\colon V(\Gamma)\rightarrow\{+,-\}\}$. Because of Lemma \ref{lem:q-isom} it is now enough to prove that the inclusion $(V'_\Gamma,0)\hookrightarrow\left(V_\Gamma,\delta_\Sigma^0\right)$ is a quasi-isomorphism.

The complex $V_\Gamma$, containing graphs with different decorating trees on a fixed core, can be described as
$$
\left(V_\Gamma,\delta^0_\Sigma\right) \simeq \left(\bigotimes_{v\in V_{\Gamma'}}(W,\delta)\right)^{G_{\Gamma'}}
$$
where $\Gamma'$ is the representative graph of the isomorphism class, $G_{\Gamma'}$ is the isomorphism group of permutations, and $(W,\delta)$ is the complex of all possible decorating trees.
Its elements are haired trees similar to graphs from $\mB_{-1}\mH_{\ge 0}\GC$ but with one special vertex called the root. Note that in $(W,\delta)$ the root can become univalent after splitting. An example of a graph from $W$ is
$$
\begin{tikzpicture}[baseline=-1ex,scale=.6]
\node[int] (a) at (0,0) {};
\node[int] (b) at (90:1) {};
\node[ext] (c) at (210:1) {};
\node[int] (d) at (330:1) {};
\node[int] (e) at (390:1) {};
\node[int] (f) at (390:2) {};
\node[int] (g) at (360:1.75) {};
\draw (b) edge (0.2,1.3);
\draw (b) edge (-0.2,1.3);
\draw (a) edge (b);
\draw (a) edge (c);
\draw (a) edge (d);
\draw (a) edge (e);
\draw (e) edge (f);
\draw (e) edge (g);
\draw (g) edge (360:2.05);
\end{tikzpicture}
$$

We do the same on $V'_\Gamma$ and get simply
$$
(V'_\Gamma,0) \simeq \left(\bigotimes_{v\in V_{\Gamma'}}(\langle\{+,-\}\rangle,0)\right)^{G_{\Gamma'}}
$$
where $\langle\{+,-\}\rangle\subseteq W$ is the 2-dimensional sub-complex spanned just by decorations $+$ and $-$, i.e.\ by
$$
e^{\pm\chi_x}(\circ)
$$
where $\circ$ is a tree that only contains the root named $x$. Now it is enough to prove that the inclusion $(\langle\{+,-\}\rangle,0)\hookrightarrow (W,\delta)$ is a quasi-isomorphism, i.e.\ that $\coH(W,\delta)$ is spanned by $+$ and $-$.

For $\Gamma\in W$ with 1-valent root let the \emph{stem} of $\Gamma$ be the sequence of vertices $x_0$ to $x_k$ such that $x_0$ is the root, $x_i$ and $x_{i+1}$ are connected by an edge, all $x_1$ to $x_{k-1}$ are 2-valent and $x_k$ is either 1-valent, 3-valent or hairy 2-valent. If the root is not 1-valent we consider the stem to be only the root.

Consider the filtration on $(W,\delta)$ on the number of edges that are not in the stem plus the number of hairs. It can be checked that the spectral sequence is regular according to \cite[Definition 5.2.10]{Wei} (c.f.\ \cite[Proposition 19]{DGC1}), and therefore converges correctly.

The first differential $\delta_0$ is the one that extends the stem and does not add new hairs. We split the complex as
$$
(W,\delta_0) = (W_0,\delta_0)\oplus(W_1,\delta_0)\oplus(W_{\geq 2},\delta)
$$
where $W_0$ is spanned by graphs with isolated root and graphs whose stem ends with hairless 1-valent vertex, $W_1$ is spanned by hairy 1-valent root and graphs whose stem ends with hairy 2-valent vertex, and $W_{\geq 2}$ is spanned by other graphs. One easily sees that $(W_{\geq 2},\delta_0)$ is acyclic, $\coH(W_0,\delta_0)$ is spanned by the isolated root and $\coH(W_1,\delta_0)$ is spanned by the hairy 1-valent root.

The two classes does not cancel on further pages, so $\coH(W,\delta)$ must be 2-dimensional. Since we already know that $+$ and $-$ represent classes, they must be the two.
\end{proof}

\begin{prop}
\label{prop:2}
    The inclusion
    $$
    \GC^+\oplus\GC^-\hookrightarrow\GC^\pm
    $$
    is a quasi-isomorphism.
\end{prop}
\begin{proof}

The complex $\GC^\pm$ splits into loops ($\mB_{0}\GC^\pm$) and non-loops ($\mB_{\ge1}\GC^\pm$), and we prove the quasi-isomorphism property for each one separately.

For the non-loops, similarly to the previous proposition, we once again introduce a structure that simplifies the way graphs are decorated while keeping the homological properties the same. We define the \emph{skeleton} of the graph as the graph obtained by recursively replacing a bivalent vertex with an edge between its neighbors. Then graphs from $\mB_{\ge1}\GC^\pm$ can be viewed as skeletons whose edges are decorated by a string of $+$ and $-$ symbols corresponding to the decorations of the original vertices (in their original order).

    We consider a filtration by the number of vertices in the skeleton. The first differential on the filtered complex is $\delta_0$ being the part of $\delta_\Sigma|_{\GC^\pm}$ that does not change the skeleton. Therefore, the associated graded complex splits as
$$
(\mB_{\ge1}\GC^\pm,\delta_0)\simeq \prod_{(\Gamma,s)}(C_{(\Gamma, s)},\delta_0)
$$
where for fixed skeleton $\Gamma$ with decorations $s$, $C_{(\Gamma,s)}$ is the complex spanned by all graphs with that skeleton and decorations of skeleton vertices.

 Similarly to the previous proposition, the complex develops into a tensor product
 $$
 C_{(\Gamma,s)}=\left(\bigotimes_{i=1}^e E_i\right)^{G_\Gamma}
 $$
 where $e$ runs through all edges in the skeleton and $E_i$ represents the chain complex of possible edge decorations. The complex $E_i$ is equivalent to the complex of finite sequences of $\{+,-\}$ with at least two elements, whose first and last element match the decorations of skeleton vertices at the ends of the skeleton edge.
 
 The differential, as induced from vertex-splitting, returns a linear combination of sequences doubling a symbol in all possible ways. Each component has a sign matching the parity of its position (induced by indexation of odd-degree edges), making the ones induced by doubling same consecutive cancel out. Therefore, for the even-length subsequences of same symbols no component survives, while for the odd-length ones only one does (the sign of which no longer matters). That property implies that the complex (denote as $A_+$) of such sequences containing only $+$ is acyclic, unless the single-element sequence $(+)$ is excluded like in $E_i$ (the same goes for $A_-$). Every other sub-complex defined by a starting symbol and the number of symbol changes (preserved by the differential), can be written as tensor product of multiple complexes $A_+$ and $A_-$ (with a single element included), making them all acyclic. Hence $E_i$ is acyclic if the first and last elements don't match, and if they do, the homology is generated by a class corresponding to $(+,+)$ or $(-,-)$.

The complex $(C_{(\Gamma, s)},\delta_0)$, as a tensor product, has trivial homology unless all the skeleton vertices have matching decorations. Therefore $\coH(\mB_{\ge1}\GC^\pm,\delta_0)$ is spanned by skeletons with all vertices having the same decorations, either $+$ or $-$. By the same argument, the same can trivially be shown for $\mB_{\ge1}\GC^+\oplus\mB_{\ge1}\GC^-$. Then Lemma \ref{lem:q-isom} implies the result (for non-loops).

As for the sub-complex of loops, it follows from our reasoning about the cohomology of $E_i$. The sub-complex $\mB_0\GC^\pm$ equivalently splits into complexes defined by the number of times the sign changes, which can then be described as the rotation co-invariants of the tensor products of $A_+$ and $A_-$ complexes. The co-invariants are acyclic because the main space is so, as we have shown. Therefore, the complement sub-complex to $\mB_0\GC^+ \oplus \mB_0\GC^-$ is acyclic, implying the result for loops.
\end{proof}

Note that the morphism $e^\chi-e^{-\chi}\colon(\GC,\delta)\rightarrow(\HoGC,\delta_\Sigma)$ can be co-restricted to $(\GC,\delta)\rightarrow(\HGC,\delta_\Sigma)$. That is the (almost) quasi-isomorphism from \cite{TW1}.

\section{Extra differentials}
\label{s:extra}

On $\fGCc$ we have an extra differential $\nabla$ that adds an edge in all possible ways. To precisely define the sigh we set that the new edge gets the last label in a representative in the space of co-invariants. In \cite[Corollary 4]{DGC1} it is shown that the complex $(\fGCc,\delta+\nabla)$ is acyclic. Since $\nabla$ cannot reduce the valence of vertices, $(\GC,\delta+\nabla)$ is a sub-complex. From the quasi-isomorphic inclusion $(\GC,\delta)\hookrightarrow(\fGCc,\delta)$ an easy spectral sequence argument implies that $(\GC,\delta+\nabla)$ is also acyclic.

In the other hand, on $\HoGC$ we have an extra differential $\Delta$ that deletes one hair and transforms it into an edge towards another vertex in all possible ways. Again set that the new edge gets the last label in a representative in the space of co-invariants. In \cite[Theorem 1.1]{DGC3} it is shown that the quotient $(\HGC,\delta+\Delta)$ is $1$-dimensional, the class is represented by $\Sigma_3$. Therefore the extended inclusion
\begin{equation}
    \K\otimes(\fGCc[-1],\delta)\hookrightarrow(\HoGC,\delta+\Delta)
\end{equation}
that sends $(\lambda,\Gamma)\mapsto\lambda\Sigma_3+\Gamma$ is a quasi-isomorphism. We extend it to quasi-isomorphisms
\begin{equation}
\label{eq:DGC3}
\K\otimes(\GC[-1],\delta)\hookrightarrow\K\otimes(\fGCc[-1],\delta)\hookrightarrow(\HoGC,\delta+\Delta).
\end{equation}

\begin{lem}
    For both signs $\pm$ it holds that:
    \begin{equation}
        e^{\pm\chi}\circ\nabla = \pm\Delta\circ e^{\pm\chi}.
    \end{equation}
\end{lem}
\begin{proof}
    Let $\Gamma\in\fGCc$. Then
    \begin{equation*}
        \begin{split}
            \Delta(e^{\pm\chi}(\Gamma))&= \sum_{i=0}^\infty\frac{(\pm 1)^i}{i!}\Delta(\chi^i(\Gamma)) \\
            &= \sum_{i=0}^\infty\frac{(\pm 1)^i}{i!}i\chi^{(i-1)}(\nabla(\Gamma))) \\
            &= \pm\sum_{i=1}^\infty\frac{(\pm 1)^{(i-1)}}{(i-1)!}\chi^{i-1}(\nabla(\Gamma))) \\
            &= \pm e^{\pm\chi}(\nabla(\Gamma)).
        \end{split}
    \end{equation*}
\end{proof}

The lemma implies that there are maps of complexes
\begin{equation}
    e^{\pm\chi}\colon(\fGCc,\delta\pm\nabla)\rightarrow(\HoGC,\delta_\Sigma+\Delta).
\end{equation}

\begin{cor}
\label{cor:main}
The cohomology $\coH(\HoGC,\delta_\Sigma+\Delta)$ is $1$-dimensional, the class being represented by $\Sigma_3$.
\end{cor}
\begin{proof}
Let us consider the map
\begin{equation}
\label{eq:map2}
    \K\oplus(\GC,\delta+\nabla)\oplus(\GC,\delta-\nabla)\rightarrow(\HoGC,\delta_\Sigma+\Delta)
\end{equation}
defined as in Theorem \ref{thm:qi}. On its mapping cone we set up a spectral sequence on the edge excess. In the first page we have a mapping cone of the map $\K\oplus(\GC,\delta)\oplus(\GC,\delta)\rightarrow(\HoGC,\delta_\Sigma)$, exactly the one shown to be quasi-isomorphism in Theorem \ref{thm:qi}. Therefore the cohomology in the first page is zero.

In the fixed degree $d$ the number of edges is fixed to $e=d-1$, and the number of vertices in a connected graph is bound by $v\leq e+1$. So the edge excess $e-v$ is also bound in all participating complexes and the spectral sequence converges correctly. Therefore the map \eqref{eq:map2} is a quasi-isomorphism.

Note that $(\GC,\delta-\nabla)$ is isomorphic to $(\GC,\delta+\nabla)$, the isomorphism being multiplication by $(-1)^{e-v}$. By the discussion following \cite[Corollary 4]{DGC1} both are acyclic. So the only remaining cohomology class in $(\HoGC,\delta_\Sigma+\Delta)\simeq\K\oplus(\GC,\delta+\nabla)\oplus(\GC,\delta-\nabla)$ comes from $\K$, being exactly $[\Sigma_3]$.
\end{proof}

\begin{thm}
\label{thm:main}
    There is a spectral sequence converging to
    $$
    \coH(\mB_{\geq 0}\HoGC[1],\delta_\Sigma+\Delta)=0
    $$
    whose $E^1$ term is
    $$
    \coH(\mB_{\geq 0}\HoGC[1],\delta+\Delta)=\coH(\GC,\delta).
    $$
    Furthermore, all differentials on odd pages are $0$, i.e.\ there are no maps between classes of $\coH(\GC,\delta)$ with opposite parity of edge excess.
\end{thm}
\begin{proof}
     Complex $(\mB_{\geq 0}\HoGC[1],\delta_\Sigma+\Delta)$ is a sub-complex of $(\HoGC[1],\delta_\Sigma+\Delta)$, and since the only cohomology class of the latter lives in the quotient $(\mB_{-1}\HoGC[1],\delta_\Sigma)$ (Corollary \ref{cor:main}), it is acyclic.

     On it we set up a spectral sequence on the edge-hair excess $e+h-v$.
     The differential on the page zero is $\delta+\Delta$, with the remaining $[\Sigma,\cdot]$ changing the edge-hair excess.

     So on the $E^1$ page we have the cohomology $\coH(\mB_{\geq 0}\HoGC[1],\delta+\Delta)$. Quasi-isomorphism \eqref{eq:DGC3} can be restricted to quasi-isomorphism $(\GC,\delta)\hookrightarrow\mB_{\geq 0}(\HoGC[1],\delta+\Delta)$, showing that on the $E^1$ page we indeed have the cohomology $\coH(\GC,\delta)$. Note that edge excess $b$ of $\coH(\GC,\delta)$ is by this quasi-isomorphism sent to edge-hair excess $b$ in $\coH(\mB_{\geq 0}\HoGC[1],\delta+\Delta)$.
     
     It can be checked that the spectral sequence is regular according to \cite[Definition 5.2.10]{Wei} (c.f.\ \cite[Proposition 19]{DGC1}), and therefore converges correctly to zero.

     Note that $[\Sigma,\cdot]$ changes edge-hair excess by even amount, so the whole complex $(\mB_{\geq 0}\HoGC[1],\delta_\Sigma+\Delta)$ is actually a direct sum of two complexes, one for even and one for odd edge-hair excess. Therefore the whole spectral sequence is the direct sum of two spectral sequences, so there can not be non-zero maps on its even pages.
\end{proof}

\begin{appendices}

\section{The picture of odd cancellations}
\label{a:pict}

In Table \ref{tbl:odd} we summarize dimensions of $\coH(\GC,\delta)$ up to degree $30$, split by edge excess. The results are taken from \cite[Table 1]{Will11Loop}, while some originally come from various earlier sources as discussed there. Green line in edge access $0$ are loop graphs, and diagonal green entries form Grothendieck-Teichm\"uller Lie algebra as shown in \cite{Will0}. Red area is known to be zero, while entries in white area is either calculated by a computer or inferred theoretically. In green area we have knowledge in further degrees, but in white area all are question marks for degrees greater that $30$.

\begin{table}[h]
{\ssmall $$
\begin{tikzpicture}[every node/.style={align=center,anchor=base,text width=2.35ex,nodes={fill=red!30}}]
\matrix (mag) [matrix of nodes,ampersand replacement=\&,row 1/.style={nodes={fill=white}},row 3/.style={nodes={fill=green!30}},column 1/.style={nodes={text width=2.35ex,fill=white}}]
{
{}\& 4 \& 5 \& 6 \& 7 \& 8 \& 9 \&10 \&11 \&12 \&13 \&14 \&15 \&16 \&17 \&18 \&19 \&20 \&21 \&22 \&23 \&24 \&25 \&26 \&27 \&28 \&29 \&30\\
-1\& 0 \& 0 \& 0 \& 0 \& 0 \& 0 \& 0 \& 0 \& 0 \& 0 \& 0 \& 0 \& 0 \& 0 \& 0 \& 0 \& 0 \& 0 \& 0 \& 0 \& 0 \& 0 \& 0 \& 0 \& 0 \& 0 \& 0 \\
|[fill=white]| 0 \& 0 \& 1 \& 0 \& 0 \& 0 \& 1 \& 0 \& 0 \& 0 \& 1 \& 0 \& 0 \& 0 \& 0 \& 0 \& 0 \& 0 \& 1 \& 0 \& 0 \& 0 \& 1 \& 0 \& 0 \& 0 \& 1 \& 0 \\
1 \& 0 \& 0 \& 0 \& 0 \& 0 \& 0 \& 0 \& 0 \& 0 \& 0 \& 0 \& 0 \& 0 \& 0 \& 0 \& 0 \& 0 \& 0 \& 0 \& 0 \& 0 \& 0 \& 0 \& 0 \& 0 \& 0 \& 0 \\
2 \& 0 \& 0 \& |[fill=green!30]| 1 \& 0 \& 0 \& 0 \& 0 \& 0 \& 0 \& 0 \& 0 \& 0 \& 0 \& 0 \& 0 \& 0 \& 0 \& 0 \& 0 \& 0 \& 0 \& 0 \& 0 \& 0 \& 0 \& 0 \& 0 \\
3 \& 0 \& 0 \& 0 \& 0 \& |[fill=green!30]| 0 \& |[fill=white]| 0 \& 0 \& 0 \& 0 \& 0 \& 0 \& 0 \& 0 \& 0 \& 0 \& 0 \& 0 \& 0 \& 0 \& 0 \& 0 \& 0 \& 0 \& 0 \& 0 \& 0 \& 0 \\
4 \&   \& 0 \& 0 \& 0 \& 0 \& 0 \& |[fill=green!30]| 1 \& |[fill=white]| 0 \& |[fill=white]| 0 \& 0 \& 0 \& 0 \& 0 \& 0 \& 0 \& 0 \& 0 \& 0 \& 0 \& 0 \& 0 \& 0 \& 0 \& 0 \& 0 \& 0 \& 0 \\
5 \&   \&   \& 0 \& 0 \& 0 \& 0 \& 0 \& 0 \& |[fill=green!30]| 0 \& |[fill=white]| 0 \& |[fill=white]| 0 \& |[fill=white]| 1 \& 0 \& 0 \& 0 \& 0 \& 0 \& 0 \& 0 \& 0 \& 0 \& 0 \& 0 \& 0 \& 0 \& 0 \& 0 \\
6 \&   \&   \&   \& 0 \& 0 \& 0 \& 0 \& 0 \& 0 \& 0 \& |[fill=green!30]| 1 \& |[fill=white]| 0 \& |[fill=white]| 0 \& |[fill=white]| 0 \& |[fill=white]| 0 \& 0 \& 0 \& 0 \& 0 \& 0 \& 0 \& 0 \& 0 \& 0 \& 0 \& 0 \& 0 \\
7 \&   \&   \&   \&   \& 0 \& 0 \& 0 \& 0 \& 0 \& 0 \& 0 \& 0 \& |[fill=green!30]| 1 \& |[fill=white]| 0 \& |[fill=white]| 0 \& |[fill=white]| 1 \& |[fill=white]| 0 \& |[fill=white]| 0 \& 0 \& 0 \& 0 \& 0 \& 0 \& 0 \& 0 \& 0 \& 0 \\
8 \&   \&   \&   \&   \&   \& 0 \& 0 \& 0 \& 0 \& 0 \& 0 \& 0 \& 0 \& 0 \& |[fill=green!30]| 1 \& |[fill=white]| 0 \& |[fill=white]| 0 \& |[fill=white]| 1 \& |[fill=white]| 0 \& |[fill=white]| 0 \& |[fill=white]| 0 \& 0 \& 0 \& 0 \& 0 \& 0 \& 0 \\
9 \&   \&   \&   \&   \&   \&   \& 0 \& 0 \& 0 \& 0 \& 0 \& 0 \& 0 \& 0 \& 0 \& 0 \& |[fill=green!30]| 1 \& |[fill=white]| 0 \& |[fill=white]| 0 \& |[fill=white]| 2 \& |[fill=white]| 0 \& |[fill=white]| 0 \& |[fill=white]| 0 \& |[fill=white]| 1 \& 0 \& 0 \& 0 \\
10\&   \&   \&   \&   \&   \&   \&   \& 0 \& 0 \& 0 \& 0 \& 0 \& 0 \& 0 \& 0 \& 0 \& 0 \& 0 \& |[fill=green!30]| 2 \& |[fill=white]| 0 \& |[fill=white]| 0 \& |[fill=white]| 2 \& |[fill=white]| 0 \& |[fill=white]| 0\& |[fill=white]| 1 \& |[fill=white]| 0 \& |[fill=white]| 0 \\
11\&   \&   \&   \&   \&   \&   \&   \&   \& 0 \& 0 \& 0 \& 0 \& 0 \& 0 \& 0 \& 0 \& 0 \& 0 \& 0 \& 0 \& |[fill=green!30]| 2 \& |[fill=white]| 0 \& |[fill=white]| 0 \& |[fill=white]| 2 \& |[fill=white]| ? \& |[fill=white]| ? \& |[fill=white]| ? \\
12\&   \&   \&   \&   \&   \&   \&   \&   \&   \& 0 \& 0 \& 0 \& 0 \& 0 \& 0 \& 0 \& 0 \& 0 \& 0 \& 0 \& 0 \& 0 \& |[fill=green!30]| 3 \& |[fill=white]| 0 \& |[fill=white]| 0 \& |[fill=white]| 4 \& |[fill=white]| ? \\
13\&   \&   \&   \&   \&   \&   \&   \&   \&   \&  \& 0 \& 0 \& 0 \& 0 \& 0 \& 0 \& 0 \& 0 \& 0 \& 0 \& 0 \& 0 \& 0 \& 0 \& |[fill=green!30]| 3 \& |[fill=white]| ? \& |[fill=white]| ? \\
14\&   \&   \&   \&   \&   \&   \&   \&   \&   \&  \&  \& 0 \& 0 \& 0 \& 0 \& 0 \& 0 \& 0 \& 0 \& 0 \& 0 \& 0 \& 0 \& 0 \& 0 \& 0 \& |[fill=green!30]| 4 \\
};
\draw (mag-2-1.north west) -- (mag-1-28.south east);
\draw (mag-1-2.north west) -- (mag-17-1.south east);
\draw[-latex, thick] (mag-3-3) edge (mag-5-4);
\draw[-latex, thick] (mag-3-7) edge (mag-7-8);
\draw[-latex, thick] (mag-3-11) edge (mag-9-12);
\draw[-latex, thick] (mag-8-13) edge (mag-10-14);
\draw[-latex, thick] (mag-3-15) edge (mag-11-16);
\draw[-latex, thick] (mag-10-17) edge (mag-12-18);
\draw[-latex, thick] (mag-3-19) edge (mag-13-20);
\draw[-latex, thick] (mag-11-19) edge (mag-13-20);
\draw[-latex, thick] (mag-12-21) edge (mag-14-22);
\draw[-latex, thick] (mag-3-23) edge (mag-15-24);
\draw[-latex, thick] (mag-13-23) edge (mag-15-24);
\draw[-latex, thick] (mag-12-25) edge (mag-16-26);
\draw[-latex, thick] (mag-14-25) edge (mag-16-26);
\draw[-latex, thick] (mag-13-26) edge (mag-15-27);
\draw[-latex, thick] (mag-3-27) edge (mag-17-28);
\draw[-latex, thick] (mag-15-27) edge (mag-17-28);
\end{tikzpicture}
$$}
\caption{\label{tbl:odd} The table of dimensions of $\coH(\GC,\delta)$.
The columns represent the degree, and rows represent the edge excess.
}
\end{table}

More precisely, the theoretical inference uses the spectral sequence from \cite{DGC1} strengthened by the result from \cite{TCon} that there are no cancellations on odd pages $1$, $3$ and $5$. Absence of those cancellations is enough for the conclusions at this stage of our knowledge of odd graph co-homomology because the hight of the white area is small enough up to degree $30$. The cancellations are depicted by arrows in Table \ref{tbl:odd}.

Our result states that there are no cancellations on all odd pages. Strictly speaking our spectral sequences is constructed from another extra differential than the one from \cite{DGC1}, and we do not technically know whether they coincide. Although they probably do, we have not pursued to check it here. For most applications it is not relevant what is the underlying extra differential, but only the existence of the spectral sequence.

\section{An example of explicit canceling}
\label{a:ex}

Here we will write explicit way of the first canceling, where the loop class generated by
\begin{equation}
\label{eq:L5}
L_5:=
\begin{tikzpicture}[baseline=-.5ex,scale=.4]
\node[int] (b) at (90:1) {};
\node[int] (c) at (162:1) {};
\node[int] (d) at (234:1) {};
\node[int] (e) at (306:1) {};
\node[int] (f) at (378:1) {};
\draw (b) edge (c);
\draw (c) edge (d);
\draw (d) edge (e);
\draw (e) edge (f);
\draw (f) edge (b);
\end{tikzpicture}
\end{equation}
is canceled with Tetrahedron class generated by
\begin{equation}
T:=
\begin{tikzpicture}[baseline=-.5ex,scale=.5]
\node[int] (b) at (90:1) {};
\node[int] (c) at (210:1) {};
\node[int] (d) at (330:1) {};
\node[int] (a) at (0:0) {};
\draw (b) edge (c);
\draw (b) edge (d);
\draw (c) edge (d);
\draw (a) edge (b);
\draw (a) edge (c);
\draw (a) edge (d);
\end{tikzpicture}\,.
\end{equation}

The graph $L_5\in\GC\subset\HoGC[1]$ generates a class in $\coH\left(\HoGC[1],\delta+\Delta\right)$. On the page $2$ of the spectral sequence of Theorem \ref{thm:main} that graph is mapped to
$$
\frac{1}{3!}\left[L_5,\Sigma_3\right] = \frac{15}{3!}\,
\begin{tikzpicture}[baseline=-.5ex,scale=.4]
\node[int] (b) at (90:1) {};
\node[int] (c) at (162:1) {};
\node[int] (d) at (234:1) {};
\node[int] (e) at (306:1) {};
\node[int] (f) at (378:1) {};
\node[int] (h) at (0,2) {};
\draw (b) edge (h);
\draw (h) edge (-0.5,2.25);
\draw (h) edge (0.5,2.25);
\draw (b) edge (c);
\draw (c) edge (d);
\draw (d) edge (e);
\draw (e) edge (f);
\draw (f) edge (b);
\end{tikzpicture}\,.
$$
The spectral sequence ends here, so $[L_5,\Sigma_3]$ represents a class in $\coH\left(\mH_{\geq 0}\GC[1],\delta+\Delta\right)$. However, $[L_5,\Sigma_3]$ has hairs, so it is not in the image of the inclusion $(\GC,\delta)\hookrightarrow\left(\mH_{\geq 0}\GC[1],\delta+\Delta\right)$.
To find what is a representative of this class in $\coH(\GC,\delta)$, we need to find another representative of the class $[\chi(L_3)]$ that is hairless. Indeed, we have
$$
(\delta+\Delta)\left(
\begin{tikzpicture}[baseline=-.5ex,scale=.4]
\node[int] (b) at (90:1) {};
\node[int] (c) at (180:1) {};
\node[int] (d) at (270:1) {};
\node[int] (e) at (0:1) {};
\node[int] (h) at (0,2) {};
\draw (b) edge (h);
\draw (h) edge (-0.5,2.25);
\draw (h) edge (0.5,2.25);
\draw (b) edge (c);
\draw (c) edge (d);
\draw (d) edge (e);
\draw (e) edge (b);
\end{tikzpicture}
\,-4\,
\begin{tikzpicture}[baseline=-.5ex,scale=.4]
\node[int] (b) at (90:1) {};
\node[int] (c) at (180:1) {};
\node[int] (d) at (270:1) {};
\node[int] (e) at (0:1) {};
\draw (b) edge (0,1.5);
\draw (b) edge (c);
\draw (c) edge (d);
\draw (d) edge (e);
\draw (e) edge (b);
\draw (c) edge (e);
\end{tikzpicture}
\right)
=
\begin{tikzpicture}[baseline=-.5ex,scale=.4]
\node[int] (b) at (90:1) {};
\node[int] (c) at (162:1) {};
\node[int] (d) at (234:1) {};
\node[int] (e) at (306:1) {};
\node[int] (f) at (378:1) {};
\node[int] (h) at (0,2) {};
\draw (b) edge (h);
\draw (h) edge (-0.5,2.25);
\draw (h) edge (0.5,2.25);
\draw (b) edge (c);
\draw (c) edge (d);
\draw (d) edge (e);
\draw (e) edge (f);
\draw (f) edge (b);
\end{tikzpicture}
-4
\begin{tikzpicture}[baseline=-.5ex,scale=.5]
\node[int] (b) at (90:1) {};
\node[int] (c) at (210:1) {};
\node[int] (d) at (330:1) {};
\node[int] (a) at (0:0) {};
\draw (b) edge (c);
\draw (b) edge (d);
\draw (c) edge (d);
\draw (a) edge (b);
\draw (a) edge (c);
\draw (a) edge (d);
\end{tikzpicture}\,,
$$
so $[L_5,\Sigma_3]\sim 60T$.
Therefore, the class $[L_5]$ cancels the class $10[T]$.\footnote{Precise signs actually depend on the ordering of edges in a representative, so pluses and minuses in this calculations are just for the illustration.}
\end{appendices}

{
\printbibliography
}

\end{document}